\documentclass[11pt,reqno,tbtags]{amsart}
\usepackage{amssymb}
\usepackage{natbib}
\bibpunct[, ]{[}{]}{;}{n}{,}{,}

\numberwithin{equation}{section}
\allowdisplaybreaks

\newtheorem*{property*}{Property \csname @currentlabel\endcsname}
\newtheorem{theorem}{Theorem}[section]

\newtheorem{corollary}[theorem]{Corollary}
\newtheorem{urn}{Urn}
\renewcommand{\theurn}{\Roman{urn}}

\theoremstyle{definition}
\newtheorem{example}[theorem]{Example}

\newtheorem{remark}[theorem]{Remark}

\theoremstyle{remark}

\newenvironment{romenumerate}{\begin{enumerate}% gives (i), (ii) etc.
 }{\end{enumerate}}

\newcounter{oldenumi}
{\setcounter{oldenumi}{\value{enumi}}
\begin{romenumerate} \setcounter{enumi}{\value{oldenumi}}}
{\end{romenumerate}}

\newcounter{thmenumerate}

\newcounter{xenumerate}   %no left indentation; thus wider lines

\newcommand{\refT}[1]{Theorem~\ref{#1}}
\newcommand{\refC}[1]{Corollary~\ref{#1}}

\newcommand{\refR}[1]{Remark~\ref{#1}}
\newcommand{\refS}[1]{Section~\ref{#1}}
\newcommand{\refU}[1]{Urn~\ref{#1}}

\newcommand\marginal[1]{\marginpar{\raggedright\parindent=0pt\tiny #1}}

\begingroup
  \divide\count255 by 60
  \multiply\count255 by -60
  \advance\count255 by \time
  \ifnum \count255 < 10 \xdef\klockan{\the\count1.0\the\count255}
  \else\xdef\klockan{\the\count1.\the\count255}\fi
\endgroup

\newcommand\set[1]{\ensuremath{\{#1\}}}

\newcommand\lrpar[1]{\left(#1\right)}

\def\rompar(#1){\textup(#1\textup)}    % usage: \rompar(...)

\def\xexp(#1){e^{#1}}

\newcommand\ntoo{\ensuremath{{n\to\infty}}}

\newcommand\ie{i.e.\spacefactor=1000}
\newcommand\eg{e.g.\spacefactor=1000}

\newcommand{\tend}{\longrightarrow}
\newcommand\dto{\overset{\mathrm{d}}{\tend}}

\newcommand\eqd{\overset{\mathrm{d}}{=}}

\newcommand\bbZ{\mathbb Z}

\newcounter{CC}
\newcounter{cc}
\renewcommand\Re{\operatorname{Re}}

\newcommand\E{\operatorname{\mathbb E{}}}
\renewcommand\P{\operatorname{\mathbb P{}}}
\newcommand\Var{\operatorname{Var}}
\newcommand\Cov{\operatorname{Cov}}

\newcommand\gb{\beta}
\newcommand\gd{\delta}

\newcommand\gl{\lambda}

\newcommand\gS{\Sigma}

\newcommand\cQ{\mathcal Q}

\newcommand\cT{{\mathcal T}}

\def\[#1]{[\![#1]\!]}

\newcommand\matrixx[1]{\begin{pmatrix}#1\end{pmatrix}}

\newcommand\qqw{^{-1/2}}
\newcommand\qw{^{-1}}

\renewcommand{\=}{:=}

\newcommand{\prt}{plane recursive tree}
\newcommand{\rprt}{random plane recursive tree}
\newcommand{\Sp}{Stirling permutation}
\newcommand{\rsp}{random Stirling permutation}
\newcommand{\gpu}{generalized \Polya{} urn}
\newcommand{\een}{second-order Eulerian number}
\newcommand{\Eul}[2]{\bigl\langle\genfrac{<}{>}{0pt}{}{#1}{#2}\bigr\rangle}
\newcommand{\tw}{trapezoidal word}
\newcommand{\dnk}{D_{nk}}
\newcommand{\snk}{S^{(k)}_n}
\newcommand{\Polya}{P\'olya}
\newcommand\REM[1]{{\raggedright\texttt{[#1]}\par\marginal{XXX}}}

\newcommand\urladdrx[1]{{\urladdr{\def~{{\tiny$\sim$}}#1}}}

\begin{document}
\title%[]
{Plane recursive trees, Stirling permutations and an urn model}

\date{March 7, 2008} % (typeset \today{} \klockan)} %; revised ...

\author{Svante Janson}
\address{Department of Mathematics, Uppsala University, PO Box 480,
SE-751~06 Uppsala, Sweden}
\email{svante.janson@math.uu.se}
\urladdrx{http://www.math.uu.se/~svante/}

%\keywords{<keywords>}
\subjclass[2000]{} 
%{60C05 (68P10,68W40)} %%{Primary: <subject>; Secondary: <subject>}

\begin{abstract} 
We exploit a bijection between \prt{s} and \Sp{s}; this yields the
equivalence of some results previously proven separately by different
methods for the two types of objects as well as some new results. We
also prove results on the joint distribution of the numbers of
ascents, descents and plateaux in a random \Sp. The proof uses an
interesting \gpu.
\end{abstract}

\maketitle

\section{Introduction}\label{S:intro}

A \emph{\prt} is a rooted plane (= ordered) tree obtained by starting
with the root and recursively adding leaves to the tree. We label the
root by 0 and the added vertices by $1,2,\dots$; a \prt{} with $n+1$
vertices is thus a labelled rooted plane tree (with labels
$0,\dots,n$) where the labels increase along each branch as we travel
from the root.
Since a new vertex may be joined to an existing vertex $i$ in $d_i+1$
positions, where $d_i$ is the outdegree of $i$, the total number of
positions for vertex $n$ is $\sum_{i=0}^{n-1}(d_i+1)=2n-1$ (and in particular
not depending on the present shape of the tree), and thus the number
of \prt{s} with $n+1$ vertices is $(2n-1)!!\=1\cdot3\dotsm(2n-1)$.
A \emph{\rprt} with $n+1$ vertices is obtained by randomly (and
uniformly) choosing one of these $(2n-1)!!$ 
trees;
equivalently, it is obtained by recursively adding $n$ vertices to the
root, each time choosing the position uniformly at random among the
possibilities. 
Random \prt{s} were studied by \citet{MSSz}.
They obtained, among other results, an asymptotic joint normal
distribution for the numbers of vertices of outdegrees 0,1,2; this was
later extended to arbitrary outdegrees by \citet{SJ155}.
More recently, the distribution of the degree of the node with a given
label has been studied by \citet{KubaP}.
See also \cite{Fla:inc} and the survey \cite{SmytheM}, where also 
other related types of trees are studied.

The well-known \emph{depth first walk} of a rooted plane tree starts
at the root, goes first to the leftmost daughter of the root, explores
that branch (recursively, using the same rules), returns to the root,
and continues with the next daughter of the root, until there are no more
daughters left. Note that every edge is passed twice in this walk, once
in each direction.
We find it convenient to label the edges too in a \prt, using the labels
$1,2,\dots$ in the order the edges are added to the tree. Thus, edge
$j$ is the edge connecting vertex $j$ to some earlier vertex.
We code a \prt{} by the sequence of the labels of the edges passed by
the depth first walk; a \prt{} with $n+1$ vertices is thus coded by a 
string of length $2n$, where each of the labels $1,\dots,n$ appears
twice. In other words, the code is a permutation of the multiset
\set{1,1,2,2,\dots,n,n}. 
Adding a new vertex $n+1$ means inserting the pair $(n+1)(n+1)$
somewhere in the code, at one of $2n+1$ possible places (first, last,
or in any of the $2n-1$ gaps between of consecutive labels)
corresponding to the $2n+1$ places in the tree where a new leaf might
be added. 
(We will in the sequel call all these $2n+1$ places 'gaps', including the
places before the first element and after the last.)
This shows that the code determines the \prt{} uniquely.
Moreover, it follows that the possible codes are exactly the
\emph{Stirling permutations} defined by \citet{GesselS}: a Stirling
permutation is a permutation of \set{1,1,2,2,\dots,n,n} such that for
each $i\le n$, the elements occuring between the two occurences of $i$
are larger than $i$. (The name comes from relations with the Stirling
numbers, see \cite{GesselS}.)

Letting $\cT_n$ be the set of \prt{s} with $n+1$ vertices (and thus $n$
edges) and $\cQ_n$ the set of Stirling permutations of length $2n$,
there is thus a simple bijection $\cT_n\cong \cQ_n$. 
We let $T_n$ denote a \rprt{} with $n+1$ vertices, \ie, a
(uniformly) random element of $\cT_n$,
and let similarly $Q_n$ denote a
random Stirling permutation of length $2n$, \ie, a (uniformly) random
element of $\cQ_n$; the bijection $\cT_n\cong \cQ_n$
thus yields a correspondence between the random objects $T_n$ and $Q_n$.
One of the purposes of this note is to show how this correspondence
connects some previous results on \rprt{s} and \rsp{s}.
We will also extend some of these results.
(The correspondence is very simple, but we have not seen it utilized
in the literature before.)

\citet{GesselS} studied the number of {descents} in \Sp{s}.
This was recently continued and extended by \citet{Bona}, using the
following definitions:
If $a_1a_2\dotsm a_{2n}$ is a \Sp, say that an index $i=0,\dots,2n$ 
(or the gap $i,i+1$)
is
an \emph{ascent} if $a_i<a_{i+1}$,
a \emph{descent} if $a_i>a_{i+1}$,
and a \emph{plateau} if $a_i=a_{i+1}$,
where we set $a_0=a_{2n+1}=0$. (Thus $0$ is always an ascent and $2n$
a descent.)
Let $X_n$, $Y_n$ and $Z_n$ denote the numbers of ascents, descents and
plateaux, respectively, in a random \Sp{} $Q_n$. Thus
\begin{equation}\label{xyz}
X_n+Y_n+Z_n=2n+1.  
\end{equation}

Let $N_{nd}$ denote the number of vertices 
with outdegree $d$
in the \rprt{} $T_n$ and let $L_n\=N_{n0}$ be the number of leaves.
It is immediately seen that in the 
correspondence above, leaves in the \prt{} correspond to plateaux in
the \Sp{}, and thus the number of leaves in the \rprt{} 
$T_n$ equals the number of plateaux in $Q_n$, \ie, 
\begin{equation}\label{lz}
 L_{n}=Z_n. 
\end{equation}

As said above, $L_n$ 
was studied by \citet{MSSz}
(note that our $L_n$ is their $L_{n+1}$). They proved, using simple recursion
relations, 
\begin{align}
  \E L_n& = \frac{2n+1}{3},
\label{EL}
\\
\Var L_n &= \frac{2(n^2-1)}{9(2n-1)}
\sim \frac{n}{9}
\label{varL}
\end{align}
and, using a generalized \Polya{} urn (see  \refS{Surn} and \refR{R2})
the asymptotic normality
\begin{equation}\label{Lnormal}
     \frac{L_n-2n/3}{\sqrt{n}}
\dto
N(0,1/9).
\end{equation}
(All unspecified limits in this paper are as \ntoo.)
Of course, in \eqref{Lnormal}, we can replace $2n/3$ by the exact mean
$(2n+1)/3$ from \eqref{EL}, and by \eqref{varL} the result can also be
written $(L_n-\E L_n)/\sqrt{\Var L_n}\dto N(0,1)$.

\citet{Bona} proved, among other things, that the random variables
$X_n$, $Y_n$ and $Z_n$ have the same distribution (reproved as
\refC{C1} below), 
and thus by \eqref{xyz}
\begin{equation}
  \label{Exyz}
\E X_n = \E Y_n = \E Z_n = \frac{2n+1}3,
\end{equation}
and further, by first showing that the probability generating
functions have only real roots,
that these random variables are asymptotically normally
distributed:
\begin{equation}
  \frac{X_n-2n/3}{\sqrt{n}}
\eqd
  \frac{Y_n-2n/3}{\sqrt{n}}
\eqd
  \frac{Z_n-2n/3}{\sqrt{n}}
\dto
N(0,1/9)
\label{Znormal}.
\end{equation}

We now see that the equality $L_n=Z_n$ means that the results 
\eqref{EL} and \eqref{Exyz} are equivalent, as well as 
\eqref{Lnormal} and \eqref{Znormal}. Furthermore, \eqref{varL} yields
also an exact formula for the variance of $X_n$, $Y_n$ and $Z_n$,
which was raised as a question in \cite{Bona}.

In \refT{T2} we will extend \eqref{Znormal} to joint convergence to
a joint normal distribution.

\begin{remark}\label{R1}
Let $C_{n,k}$ be the number of \prt{s} with $n+1$ vertices and $k$
leaves, or, equivalently by the discussion above, the number of \Sp{s}
of length $2n$ with $k$ plateaux (or $k$ ascents, or $k$ descents). It
is easy to see that $C_{n,k}$ satisfies the recursion
\begin{equation}
  \label{rec}
C_{n,k}
=kC_{n-1,k}+(2n-k)C_{n-1,k-1}
\end{equation}
for all $n\ge2$ and $k\ge1$ (or $k\in\bbZ$), 
with $C_{1,k}=\gd_{1k}$ and $C_{n,0}=0$,
see
\cite{GesselS}, \cite{MSSz}, \cite{CM}.
These numbers $C_{n,k}$ are called \een{s} \cite[\S6.2]{CM}; in
standard notation $C_{n,k}=\Eul{n}{k-1}$ 
\cite{CM},
and thus
\cite{MSSz}, for $n\ge1$,
\begin{equation}
  \P(L_n=k)
=
  \P(X_n=k)
=
  \P(Y_n=k)
=
  \P(Z_n=k)
=
\frac{\Eul{n}{k-1}}{(2n-1)!!}.
\end{equation}  
\end{remark}

\begin{remark}
Another combinatorial interpretation of the same numbers $C_{n,k}$ is
given by \citet{Riordan} as the number of \emph{\tw s} of length $n$
with $k$ distinct elements, where a \tw{} is a word $a_1\dotsm a_n$
with $a_i\in\set{1,2,\dots,2i-1}$, \ie, an element of
$[1]\times[3]\times\dotsm\times[2n-1]$. 
(Again, it is easy to verify the recursion \eqref{rec}.)
Hence $L_n$ etc.\ also give the distribution of the number of distinct
elements in a random \tw{} of length $n$.
We do not know any interesting statistics of \tw{s} that correspond to
other statistics of \prt{s} or \Sp{s} such as $N_{nd}$ ($d\ge1$) or
the triple $(X_n,Y_n,Z_n)$.
\end{remark}

\begin{remark}
For the connections between the \een{s} $C_{n,k}$ and Stirling numbers,
see \eg{}
\cite{Ginsburg},
\cite{Carlitz},
\cite{Riordan},
\cite{GesselS},
\cite{CM}.
\end{remark}

The correspondence between \prt{s} and \Sp{s} makes it natural to
study also other statistics of one of these objects and see what they
correspond to for the other object, thus giving more equalities of the
type \eqref{lz}. (In some cases the results are, however, disappointing
in that the resulting statistics seem uninteresting.)

\begin{example}
  The number of plateaux in a \Sp{} is, as discussed above, equal to
  the number of leaves in the corresponding \prt.
It seems natural to consider ascents and descents also. It is
  straight-forward to see that an ascent corresponds to a non-root vertex in
  the tree that either has no sister to its left, or else has a higher
  label than the sister that is nearest to it to the left;
similarly,
a descent corresponds to a non-root vertex in
  the tree that either has no sister to its right, or else has a higher
  label than the sister that is nearest to it to the left;
We can thus interpret the results on $(X_n,Y_n,Z_n)$ as results on the
  numbers of such vertices in a \rprt, but it remains to show that
  these numbers are interesting.
\end{example}

\begin{example}
The distance $\dnk$ between the two occurences of $k$ in a random \Sp{}
of length $2n$ is perhaps a more interesting example. It
equals $2\snk-1$ where $\snk$ is the size of the subtree rooted at
vertex $k$ of the corresponding random \prt. The latter variable was studied by
\citet{MSSz} using a \Polya--Eggenberger urn;
they obtained both exact
formulas for the distribution, mean and variance, as well as an asymptotic
beta distribution: $\snk/n\dto\gb(\tfrac12,k)$ as \ntoo{} for every
fixed $k\ge1$. 
(Recall that we label the root by 0, so the labels are shifted from
\cite{MSSz}.) 
 These results immediately tranfer to results about $\dnk$, for
 example,
for any fixed $k\ge1$,
\begin{equation}
  \frac{\dnk}{2n}
\dto
\gb(\tfrac12,k)
\qquad
\text{as \ntoo}.
\end{equation}
\end{example}

\begin{example}
  The degree of the root of $T_n$ was also studied in \cite{MSSz}.
This corresponds for a \Sp{} $Q_n$ to the number $d$ such that the 
\Sp{} is of the form
$a_1\dotsm a_1 a_2\dotsm a_2 \dotsm a_d\dotsm a_d$.
Again, this is not of any obvious great interest.
\end{example}

\section{An interesting urn model}\label{Surn}

\citet{MSSz} used a representation with a \gpu{} to prove the
asymptotic normality \eqref{Lnormal}. 
(In fact, they proved joint asymptotic normality of 
$N_{nd}$ ($d\le2$), 
which was extended in \cite{SJ155} to all $d$, see also
\cite[Example 7.6]{SJ154}.)
We use a similar urn to study the joint distribution of $(X_n,Y_n,Z_n)$.

If we extend a \Sp{} in $\cQ_n$ by inserting the pair $(n+1)(n+1)$ at
one of the $2n+1$ possible places (gaps), then the ascent, descent or plateau
at that gap is destroyed and replaced by the sequence
\emph{ascent, plateaux, descent}
at the three resulting new gaps.
Consequently, the random vector $(X_n,Y_n,Z_n)$ is described by the
following \gpu{} model:

\begin{urn}\label{UA}
Consider an urn with balls of three colours, and let $(X_n,Y_n,Z_n)$
be the number of balls of each colour at time $n$. At each time step,
draw one ball at random from the urn, discard it, and add one new ball
of each colour. Start with $(X_1,Y_1,Z_1)=(1,1,1)$.
\end{urn}

This urn model is completely symmetric in the three colours, and we
thus immediately see the following.

\begin{theorem}
  \label{T1}
For each $n\ge1$, the distribution of\/ $(X_n,Y_n,Z_n)$ is exchangeable,
\ie, invariant under any permutation of the three variables.
\end{theorem}

\begin{corollary}[\citet{Bona}]
  \label{C1}
$X_n\eqd Y_n \eqd Z_n$.
\end{corollary}

\begin{remark}
  We have stated \refT{T1} and \refC{C1} for a fixed $n$, but the
  results extend to the processes $(X_n)_n$, $(Y_n)_n$, $(Z_n)_n$.
\end{remark}

It is customary and convenient to formulate \gpu{s} using drawings
with replacement. In the case of \refU{UA}, 
we thus restate the description above and say instead
that we draw a ball and
replace it together with one ball each of the two other colours.
In other words, \refU{UA} is described by the replacement matrix
\begin{equation}\label{A}
A=  \matrixx{0 & 1 & 1 \\ 1 & 0 & 1 \\ 1 & 1 & 0}.
\end{equation}

We now easily obtain one of our main results.

\begin{theorem}\label{T2}
  $(X_n,Y_n,Z_n)$ are jointly asymptotically normal:
\begin{equation}\label{t2a}
n\qqw{(X_n-2n/3,Y_n-2n/3,Z_n-2n/3)}
\dto
N\lrpar{0,\Sigma},
\end{equation}
where the asymptotic covariance matrix is given by
\begin{equation}\label{t2b}
\Sigma= \frac1{18} 
\matrixx{\phantom{-}2 & -1 & -1 \\ 
  -1 & \phantom{-}2 & -1 \\ 
  -1 & -1 & \phantom{-}2}.
\end{equation}
\end{theorem}

\begin{proof}
It is easily seen that the matrix $A$ in \eqref{A} has eigenvalues
$2$, $-1$, $-1$. 
(Note also that $A$ can be regarded as convolution 
with $(0,1,1)$
on the
group $\bbZ_3$, and thus the eigenvalues are given by the Fourier
transform of this vector.)
In particular, the dominant eigenvalue $\gl_1=2$ (the row sum in $A$),
and all other eigenvalues lie in the half plane \set{\Re\gl<\gl_1/2}.
Consequently, \cite[Theorem 3.22]{SJ154} applies and shows the 
asymptotic normality \eqref{t2a}.

To identify the asymptotic covariance matrix $\gS$, we may use
\cite[Lemma 5.4 and (2.15), or Lemma 5.3]{SJ154} and some
  straight-forward calculations. We may, alternatively, avoid
  calculations completely by noting that the diagonal entries, which
  must be equal by \refC{C1}, are $1/9$ by 
\eqref{Znormal}; furthermore, the off-diagonal entries also are equal
by \refT{T1}, and the row sums in $\gS$ are 0 as a consequence of
\eqref{xyz}; hence the off-diagonal entries are $-1/18$.
\end{proof}

An exact formula for the covariances is easily obtained too.
\begin{theorem}
For every $n\ge1$,
\begin{equation}
  \Cov(X_n,Y_n)
=
  \Cov(X_n,Z_n)
=
  \Cov(Y_n,Z_n)
=
- \frac{(n^2-1)}{9(2n-1)}.
\end{equation}
\end{theorem}
\begin{proof}
  This can be proved using recursion formulas derived from the urn
  model.
It is, however, simpler to use \refT{T1} and the already known
  variance
\eqref{varL}. 
Indeed, by symmetry, the three covariances are equal.
Further, by \eqref{xyz},
$$\Var(X_n)+\Cov(X_n,Y_n)+\Cov(X_n,Z_n)=\Cov(X_n,X_n+Y_n+Z_n)=0.$$
Since $X_n\eqd Z_n=L_n$, we thus obtain from \eqref{varL} 
\begin{equation*}
  \Cov(X_n,Y_n)=-\frac12\Var(X_n)=
-\frac12\Var(L_n)=
- \frac{(n^2-1)}{9(2n-1)}.
\qedhere
\end{equation*}
\end{proof}

\begin{remark}\label{R2}
  If we only are interested in the univariate distribution of the
  variables,  
we may combine two of the colours into one, and thus instead
  study 
Urn
\stepcounter{urn}\theurn:
the two-colour urn with replacement vectors $(0,2)$ and
  $(1,1)$.
This is the urn used by \cite{MSSz} to show the asymptotic normality
\eqref{Lnormal}.
\end{remark}

\begin{remark}
  \refU{UA} is a rather special type of \gpu, since (using the
  formulation of drawing without replacement), the added balls do not
  depend on the drawn ball. This is perhaps even more striking in the
  corresponding continuous-time multitype branching process 
(see  \cite{AK}, \cite[\S V.9]{AN}, \cite{SJ154}), which now can be
  described as follows: 
\emph{%  
There are individuals of three types (colours). Each individual has an
exponential lifetime, independent of all other individuals. When
someone dies, one new individual of each of the three colours is born.}
Since it does not matter which individuals that die, it might be
thought that the result would be like taking a large number $3n$ of
individuals, $n$ of each colour, and randomly removing $n-1$
of them. However, this is not correct, not even asymptotically; in
fact, a simple calculation of the asymptotic variance in this
simplified model (where the number of remaining individuals of a given
colour has a hypergeometric distribution) yields $(4/27)n$ instead of $n/9$.
The reason is that although the individuals die independently, they
  are born together in triplets. Moreover, since the individuals in
  older triplets have larger probabilities of having died at some
  given instance, there is a positive correlation between the deaths
  (up to a given time) for the individuals in the same triplet; since
  these individuals have different colours, this tends to decrease the
  variances.

Consequently, we drawn the conclusion (and warning!)\ that although
\refU{UA} is a very simple type of \gpu, the fact that the replacements do
not depend on the drawn colour does not really simplify the arguments.
(At least, we do not see any simplification.) This also emphasizes
that it is usually better to formulate \gpu{} models using drawing
with replacement; in our case we than have the replacement matrix $A$
in \eqref{A},
which (although very nice in other respects) describes replacements
that do depend on the drawn colour.
\end{remark}

\section{Further comments}

The proof by \citet{Bona} of the asymptotic normality \eqref{Znormal}
is completely different and is based on first showing that the probability
generating functions have only real roots; thus the random variables
can be represented as sums of independent Bernoulli variables. 
(These Bernoulli variables have in general irrational means and do not
have any combinatorial interpretation.)
This is a strong property
that implies not only asymptotic normality (provided the variance tends to
infinity, as in this case); it has other desirable consequences too,
for example it leads to explicit error estimates for the convergence
to the normal distribution as well as to large deviation estimates
(Chernoff bounds). 

The success of the generating function in this context suggests that
it might be interesting and profitable to study the trivariate probability
generating function for $(X_n,Y_n,Z_n)$.

Note that, by \refR{R1}, the univariate 
probability generating function of $L_n$ (or $Z_n$) is
$(2n-1)!!\qw$ times the generating function $\sum_k C_{n,k}x^k$, which
has interesting properties studied in 
\cite{Ginsburg},
\cite{Carlitz},
\cite{Riordan},
\cite{GesselS},
\cite{CM}. It might be hoped that the trivariate generating function too
has some interesting combinatorial properties.

\newcommand\AAP{\emph{Adv. Appl. Probab.} }
\newcommand\JAP{\emph{J. Appl. Probab.} }
\newcommand\JAMS{\emph{J. \AMS} }
\newcommand\MAMS{\emph{Memoirs \AMS} }
\newcommand\PAMS{\emph{Proc. \AMS} }
\newcommand\TAMS{\emph{Trans. \AMS} }
\newcommand\AnnMS{\emph{Ann. Math. Statist.} }
\newcommand\AnnPr{\emph{Ann. Probab.} }
\newcommand\CPC{\emph{Combin. Probab. Comput.} }
\newcommand\JMAA{\emph{J. Math. Anal. Appl.} }
\newcommand\RSA{\emph{Random Struct. Alg.} }
\newcommand\ZW{\emph{Z. Wahrsch. Verw. Gebiete} }
\newcommand\SPA{\jour{Stochastic Process. Appl.} } 
\newcommand\DMTCS{\jour{Discr. Math. Theor. Comput. Sci.} }

\newcommand\AMS{Amer. Math. Soc.}
\newcommand\Springer{Springer-Verlag}
\newcommand\Wiley{Wiley}

\newcommand\vol{\textbf}
\newcommand\jour{\emph}
\newcommand\book{\emph}
\newcommand\inbook{\emph}
\def\no#1#2,{\unskip#2, no. #1,} %(typeset after year) 
\newcommand\toappear{\unskip, to appear}

\newcommand\webcite[1]{%\hfil  %???
  %\penalty0 %???
\texttt{\def~{{\tiny$\sim$}}#1}\hfill\hfill}
\newcommand\webcitesvante{\webcite{http://www.math.uu.se/~svante/papers/}}
\newcommand\arxiv[1]{\webcite{arXiv:#1.}}

\def\nobibitem#1\par{}


\begin{thebibliography}{99}

\bibitem{AK}
K. B. Athreya \& S. Karlin,
Embedding of urn schemes into continuous time Markov branching
processes and related limit theorems.
\AnnMS \vol{39} (1968), 1801--1817.

\bibitem{AN}
K. B. Athreya \& P. E. Ney,
\emph{Branching Processes}.
\Springer, Berlin, 1972.

\bibitem[Bergeron, Flajolet and Salvy(1992)]{Fla:inc}
F. Bergeron, P. Flajolet \& B. Salvy,
Varieties of increasing trees.
\inbook{CAAP '92 (Rennes, 1992),
Lecture Notes in Comput. Sci.} {581},
Springer, Berlin, 1992,
24--48.


\bibitem[Bona(2007+)]{Bona}
M. Bona,
Real zeros and normal distribution for statistics on Stirling
permutations defined by Gessel and Stanley. 
Preprint, 2007.
\arxiv{0708.3223v1} % [math.CO]

\bibitem[Carlitz(1965)]{Carlitz}
L. Carlitz, 
The coefficients in an asymptotic expansion.  
\emph{Proc. Amer. Math. Soc.}  \vol{16}  (1965) 248--252.


\bibitem[Gessel and Stanley(1978)]{GesselS}
I. Gessel \& R. P. Stanley, 
Stirling polynomials.  
\emph{J. Combinatorial Theory Ser. A}  \vol{24}  (1978), no. 1, 24--33.

\bibitem[Ginsburg(1928)]{Ginsburg}
J. Ginsburg, %Jekuthiel
%Questions and Discussions: Discussions: 
Note on Stirling's Numbers.
\jour{Amer. Math. Monthly} \vol{35} (1928), no. 2, 77--80. 

\bibitem[Graham, Knuth and Patashnik(1994)]{CM}
R.L. Graham, D.E. Knuth \& O. Patashnik,
\emph{Concrete Mathematics}.
2nd ed.,
Addison--Wesley, Reading, Mass., 1994.

\bibitem{SJ154}
S. Janson,
Functional limit theorems for multitype branching processes
and generalized P\'olya urns.
\SPA \vol{110} (2004), no. 2,  177--245.

\bibitem[Janson(2005)]{SJ155}
S. Janson,
Asymptotic degree distribution in random recursive trees.
\RSA \vol{26} (2005), no. 1--2, 69--83.  

\nobibitem{KnuthI} 
D.E. Knuth, 
\emph{The Art of Computer Programming. Vol. 1:
 Fundamental Algorithms}. 
3nd ed., Addison-Wesley,
Reading, Mass., 1997.

\bibitem[Kuba and Panholzer(2007)]{KubaP}
M. Kuba \& A. Panholzer,
On the degree distribution of the nodes in increasing trees.  
\emph{J. Combin. Theory Ser. A}  \vol{114}  (2007),  no. 4, 597--618.


\bibitem[Mahmoud, Smythe and Szyma\'nski(1993)]{MSSz}
H. M. Mahmoud,  R. T. Smythe \& J. Szyma\'nski,  
On the structure of random plane-oriented recursive trees and their
branches. 
\RSA \vol4 (1993), no. 2, 151--176.

\bibitem[Riordan(1976)]{Riordan}
J. Riordan, 
The blossoming of Schr\"oder's fourth problem.
\emph{Acta Math.} \vol{137} (1976), no. 1--2, 1--16. 

\bibitem[Smythe and Mahmoud(1995)]{SmytheM}
R. T. Smythe and H. Mahmoud, 
A survey of recursive trees.
\emph{Theory Probab. Math. Statist.}
\vol{51} (1995), 1--27.
\end{thebibliography}
\end{document}